\documentclass{article}     % Basic GT/GTM/AGT style
\usepackage{amssymb}
\usepackage[english]{babel}
\title{(Non)-completeness of $\R$--buildings and fixed point theorems}

\author{Koen Struyve \\ \\
Department of Pure Mathematics and Computer Algebra \\
Ghent University \\
Krijgslaan 281 - S22 \\
B-9000 Gent \\
Belgium \\
kstruyve@cage.ugent.be}

\parskip\medskipamount
\newtheorem{thm}{Theorem}[section]
\newtheorem{Theorem}{Main Result}
\newtheorem{rem}[thm]{Remark}

\newtheorem{lemma}[thm]{Lemma}
\newtheorem{cor}[thm]{Corollary}
\newtheorem{con}[thm]{Conjecture}

\def\<{\langle}
\def\>{\rangle}
\newcommand{\proof}{\emph{Proof.~}}

\newcommand{\dd}{\mathsf{d}}

\newcommand{\cF}{\mathcal{F}}

\def\qed{{\hfill\hphantom{.}\nobreak\hfill$\Box$}}

\newcommand{\A}{\mathbb{A}}
\newcommand{\R}{\mathbb{R}}
\newcommand{\N}{\mathbb{N}}

\begin{document}

\maketitle

%%%%%%%%%%%%%%%%%%%%   Start of main body of article
\begin{abstract}    % type your abstract below
We prove two generalizations of results appearing in~\cite{Bru-Tit:72} involving metrical completeness and $\R$--buildings. Firstly, we give a generalization of the Bruhat--Tits fixed point theorem also valid for non-complete $\R$--buildings, but with the added condition that the group is finitely generated.

Secondly, we generalize a criterion which reduces the problem of completeness to the wall trees of the $\R$--building. This criterion occured in~\cite{Bru-Tit:72} for $\R$--buildings arising from root group data with valuation.

\end{abstract}

\section{Introduction}
In 1986, Jacques Tits classified affine buildings of rank at least four (\cite{Tit:86}). He also included in this classification non-discrete generalizations of affine buildings, these metric spaces are called non-discrete affine apartment systems or $\R$--buildings. Although the first definition for these $\R$--buildings only appeared in the aforementioned paper, the examples that arise from root group data with a (non-discrete) valuation were already studied in 1972 in the book~\cite{Bru-Tit:72} by Bruhat and Tits. The classification of Tits shows that an $\R$--building of rank at least four (or equivalently dimension at least three) necessarily arises from such a root group datum with valuation. For the dimension two cases there exist various explicit and `free' constructions (for example see~\cite{Ber-Kap:*} and~\cite{Str-Mal:*}). The dimension one $\R$--buildings are also known as $\R$--trees.

As $\R$--buildings are metric spaces, one has the notion of (metrical) completeness. Affine buildings are always complete due to their discrete nature. For non-discrete $\R$--buildings completeness is rarer. Perhaps the easiest examples of fields with non-discrete valuation from which one can define complete $\R$--buildings are the already quite involved Hahn--Mal'cev--Neumann series. 

The current paper is in some sort a generalization of two results appearing in the book~\cite{Bru-Tit:72}, both involving completeness. The first such result is the Bruhat--Tits fixed point theorem. It says that a bounded group of isometries of a complete $\R$--building (or more generally a complete CAT(0)--space) has a fixed point. Our result proves the existence of a fixed point also in the non-complete case, but with the added condition that the group is finitely generated.

The second result is about the question when an $\R$--building is (metrically) complete and when not. In~\cite{Bru-Tit:72} it is proven that an $\R$--building arising from a root group datum with valuation is complete if and only if its corresponding wall trees (which are $\R$--trees) are complete. So this reduces the problem to the easier one dimensional case. We generalize this result to all $\R$--buildings, not only those arising from a root group datum with valuation, by giving a geometric proof instead of an algebraic proof.

\textbf{Acknowledgement.} The author is supported by  the Fund for Scientific Research --
Flanders (FWO - Vlaanderen). The author also would like to thank Linus Kramer for helpful remarks on how to embed non-complete $\R$-buildings into complete ones (see Lemma~\ref{embed}).

\section{Preliminaries}
\subsection{$\R$--buildings}
\subsubsection{Definition}
Let $(\overline{W},S)$ be a spherical Coxeter system of rank $n$. The group $\overline{W}$ can be realized as a finite reflection group acting on an $n$-dimensional real affine space $\A$, called the \emph{model space}. A \emph{wall} of $\A$ is a hyperplane of it fixed by a conjugate of an involution in $S$. A \emph{root} is a (closed) half-space of $\A$ bordered by a wall. The set of all walls of $\A$ defines a poset of simplicial cones in $\A$ (called \emph{sector-faces}), which forms the simplicial complex of the Coxeter system $(\overline{W},S)$. The maximal cones are called \emph{sectors}, the one less than maximal the \emph{sector-panels}. The apex of a cone formed by a sector-face in $\Lambda$ is called the \emph{base point} of that sector-face. (See~\cite[Chapter 1]{Abr-Bro:08} for a detailed discussion on finite reflection groups.)

Let $W$ be the group acting on $\A$ generated by $\overline{W}$ and the translations of $\A$. 

Consider a pair $(\Lambda,\cF)$ where $\Lambda$ is the set of \emph{points} forming a metric space together with a metric $\dd$, and $\cF$ a set of isometric injections (called \emph{charts}) from the model space $\A$ (equipped with the Euclidean distance) into $\Lambda$. An image of the model space is called an \emph{apartment}, an image of a root a \emph{half-apartment} and an image of a sector(-face/panel) is called again a \emph{sector(-face/panel)}. The pair $(\Lambda,\cF)$ is an \emph{$\R$--building} if the following five properties are satisfied:
\begin{itemize}
 \item[(A1)] If $w\in W$ and $f\in \cF$, then $f \circ w \in \cF$.
 \item[(A2)] If $f,f' \in \cF$, then $X=f^{-1}(f'(\A))$ is a closed and convex subset of $\A$, and $f|_X = f'\circ w|_X$ for some $w\in W$.
 \item[(A3)] Each two points of $\Lambda$ lie in a common apartment.
\end{itemize}
This last axiom implies that the metric $\dd$ on $\Lambda$ is defined implicitly by the isometric injections $\cF$.
\begin{itemize}
 \item[(A4)] Any two sectors $S_1$ and $S_2$ contain subsectors $S_1' \subset S_1$ and $S_2' \subset S_2$ lying in a common apartment.
 \item[(A5)] If three apartments intersect pairwise in half-apartments, then the intersection of all three is non-empty.
\end{itemize}

\subsubsection{Global and local structure}
Two sector-faces are \emph{parallel} if the Hausdorff distance between both is finite. This relation is an equivalence relation due to the triangle inequality. The equivalence classes (named \emph{simplices at infinity}, or the \emph{direction} $F_\infty$ of a sector-face $F$) form a spherical building $\Lambda_\infty$ of type $(\overline{W},S)$ called the \emph{building at infinity} of the $\R$--building $(\Lambda,\cF)$. The chambers of this building are the equivalence classes of parallel sectors. Each apartment $\Sigma$ of $(\Lambda,\cF)$ corresponds to an apartment $\Sigma_\infty$ of $\Lambda_\infty$ in a bijective way.

Two sector-faces are \emph{asymptotic} if there exists a subsector-face of both having the same dimension as the original two. Asymptotic sector-faces are necessarily parallel, the inverse is only true for sectors (see~\cite[Cor. 1.6]{Par:00}). Asymptoticness is an equivalence relation as well.

One can also define local equivalences. Let $\alpha$ be a point of $\Lambda$, and $F,F'$ two sector-facets based at $\alpha$. Then these two sector-facets will \emph{locally coincide} if their intersection is a neighbourhood of $\alpha$ in both $F$ and $F'$. This relation forms an equivalence relation defining \emph{germs of facets} as equivalence classes (notation
$[F]_\alpha$). These germs form a (weak) spherical building  $\Lambda_\alpha$ of type $(\overline{W},S)$,
called the \emph{residue} at $\alpha$. 

A detailed study of $\R$--buildings can be found in~\cite{Par:00}. We list some results from that paper here in order to refer to them later:

 \begin{lemma}[\cite{Par:00}, Prop. 1.8]\label{p8}
Let $C$ be a chamber of the building at infinity $\Lambda_\infty$ and $S$ a sector based at $\alpha \in
\Lambda$. Then there exists an apartment $\Sigma$ containing an element of the germ $[S]_\alpha$ and such that $\Sigma_\infty$ contains $C$. 
 \end{lemma}

 \begin{cor}[\cite{Par:00}, Cor. 1.9]
Let $\alpha$ be a point of $\Lambda$ and $F_\infty$ a facet of the building at infinity. Then there is a unique
facet $F' \in F_\infty$ based at $\alpha$. 
 \end{cor}

The unique facet of the previous corollary will be denoted by $F_\alpha$. We will often use the subtext to indicate the base point of the sector-face. An exception is when we use the symbol $\infty$ as subtext to denote the direction of the sector-face. To give an example, if one says $S_\infty$, $S_\alpha$ and $S_\beta$; then $S_\infty$ denotes some simplex at infinity and $S_\alpha$, $S_\beta$ are the unique parallel sector-faces based at respectively $\alpha$ and $\beta$ with direction $S_\infty$.

\begin{lemma}[\cite{Par:00}, Prop. 1.12]\label{p12}
If the germs $[S]_\alpha$ and $[S']_\alpha$ of two sectors $S_\alpha$ and $S'_\alpha$ form opposite chambers of the residue $\Lambda_\alpha$, then there exists a unique apartment containing both sectors.
\end{lemma}

\begin{lemma}[\cite{Par:00}, Prop. 1.17]\label{p17}
Let $\Sigma$ be an apartment and $S_\alpha$ a sector in $\Sigma$ based at some point $\alpha$. There exists a retraction $r$ of $\Lambda$ on $\Sigma$ such that $r$ preserves the distance of points of $\Lambda$ to $\alpha$ and does not increase other distances. Also each sector based at $\alpha$ is mapped isometrically to a sector in $\Sigma$. The only sector mapped to $S_\alpha$ is $S_\alpha$ itself.
\end{lemma}

%\begin{rem}\em
%In order to avoid pathological cases, we assume that the spherical building at infinity is thick. In particular this implies that the type of this spherical building cannot contain a direct factor $\mathsf{H}_3$ or $\mathsf{H}_4$. Non-thick cases can be reduced to thick cases, see~\cite[Prop 4.9.2]{Kle-Lee:97}.
%\end{rem}

\subsubsection{Wall and panel trees}
With a wall $M$ of an $\R$--building one can associate a direction at infinity (by taking the direction of all sector-facets it contains). This direction $M_\infty$ at infinity will be a wall of the spherical building at infinity. 

Let $m$ (respectively $\pi$) be a wall (resp. a panel contained in the wall $m$) of the building at infinity. Let $T(m)$ be the set of all walls $M$ of the $\R$--building with $M_\infty = m$, and $T(\pi)$ the set of all asymptotic classes of sector-panels in the parallel class of $\pi$. 

One can define charts (and so also apartments) from $\R$ to $T(m)$ (resp. $T(\pi)$). First choose $M$ (resp. $D$) a wall (resp. a sector-panel contained in $M$) of the model space such that there exists some chart $f$ such that $f(M)_\infty=m$ and $f(D) \in \pi$. One can identify the model space $\A$ with the product $\R \times M$. For every chart $g \in \cF$ of the $\R$--building $(\Lambda,\cF)$ such that $g(M)_\infty =m$ (resp. $f(D) \in \pi$), one defines a chart $g'$ as follows: if $x \in \R$, then $g'(r)$ is the wall $g(\{r\} \times M)$ (resp. the asymptotic class containing $g(\{r\} \times D)$).

These two constructions yield $\R$--buildings with a rank one building at infinity, such $\R$--buildings are better know as \emph{$\R$--trees} (or shortly trees when no confusion can arise). The following theorem shows the connection between both constructions.

\begin{thm}[\cite{Tit:86}, Prop. 4]
If $\pi$ is a panel in some wall $m$ at infinity, then for each asymptotic class $D$ of sector-panels with direction $\pi$, there is a unique wall $M$ in the direction $m$ containing an element of $D$. The corresponding map $D \mapsto M$ is an isometry from the $\R$--tree $T(\pi)$ to the $\R$--tree $T(m)$. 
\end{thm}

The trees obtained from walls (resp. panels) are called \emph{wall trees} (resp. \emph{panel trees}).

\subsection{CAT(0)--spaces}

For now suppose that $(X,\dd)$ is some metric space, not necessarily an $\R$--building. A \emph{geodesic} is a subset of the metric space $X$ isometric to a closed interval of real numbers. The metric space $(X,\dd)$ is a \emph{geodesic metric space} if each two points of $X$ can be connected by a geodesic. From property (A3) it follows that $\R$-buildings are geodesic metric spaces. 

Let $x,y$ and $z \in X$ be three points in a geodesic metric space $(X,\dd)$. Because of the triangle inequality we can find three points $\bar{x}, \bar{y}$ and $\bar{z}$ in the Euclidean plane $\R^2$ such that each pair of points has the same distance as the corresponding pair in $x,y,z$. The triangle formed by the three points is called a \emph{comparison triangle} of $x,y$ and $z$. Consider a point $a$ on a geodesic between $x$ and $y$. So we have that $\dd(x,y) = \dd(x,a) + \dd(a,y)$ (note that the geodesic, and so also the point $a$, is not necessarily unique). We now can find a point $\bar{a}$ on the line through $\bar{x}$ and $\bar{y}$ such that the pairwise distances in $\bar{x},\bar{y},\bar{a}$ are the same as in $x,y,a$. If the distance between $z$ and $a$ is smaller than the distance between $\bar{z}$ and $\bar{a}$ for all possible choices of $x,y,z$ and $a$, we say that the geodesic metric space $(X,d)$ is a \emph{CAT(0)--space}. %One can roughly think of CAT(0)--spaces as metric spaces with non-positive curvature, see~\cite{Bri-Hae} for more information.

A metric space is \emph{complete} if all Cauchy sequences converge. A group of isometries acting on a metric space is \emph{bounded} if at least one orbit (and hence all orbits) is a bounded set. Finite groups of isometries are always bounded.

The metric spaces formed by $\R$--buildings are examples of CAT(0)--spaces. For complete CAT(0)--spaces one has the following important theorem known as the Bruhat--Tits fixed point theorem.
\begin{thm}[\cite{Bru-Tit:72}, Prop. 3.2.4]\label{bt}
If $G$ is a bounded group of isometries of a complete CAT(0)--space $(X,\dd)$, then $G$ fixes some point of $X$.
\end{thm}

\begin{rem} \em
The notion of completeness has also another meaning when used for $\R$--buildings, in the sense of `the complete system of apartments'. However, there will be no confusion possible as we will not use this other notion.
\end{rem}

\subsection{Convex sets in spherical buildings}
Consider a (weak) spherical building $\Delta$ of type $(\overline{W},S)$ as a chamber complex. Between two chambers of the building one can define a $\overline{W}$--valued function $\delta$, called the \emph{Weyl distance} (see~\cite[Section 4.8]{Abr-Bro:08}). The word length with respect to the generating set $S$ makes this distance into a  A chamber subcomplex $K$ of $\Delta$ is \emph{convex} in $\Delta$ if each minimal gallery between two chambers of $K$ also lies in $K$. An equivalent definition is if $C$ and $D$ are chambers of $K$, and $E$ is a chamber of $\Delta$ such that $\delta(C,D)=\delta(C,E)\delta(E,D)$, then $E$ is a chamber of $K$ as well.

The following theorem is known as the center conjecture:
\begin{con}\label{ccj}
Let $\Delta$ be a (weak) spherical building and $K$ a convex chamber subcomplex of $\Delta$. Then (at least) one of the following two possibilities holds.
\begin{itemize}
\item
For each chamber $C$ in $K$ there is a chamber $D$ in $K$ opposite to $C$.
\item
The group of automorphisms of $\Delta$ stabilizing $K$ stabilizes a non-trivial simplex of $K$.
\end{itemize} 
\end{con} 

The center conjecture is most often stated in a more general way to include subcomplexes of lower rank than $\Delta$, but we omit this as we will not need it. Although it is called a conjecture, it has been proven except for the case where one has a direct factor of type $\mathsf{H}_4$. The cases $\mathsf{A}_n$, $\mathsf{B}_n$, $\mathsf{C}_n$ and $\mathsf{D}_n$ have been proved by Bernhard M\"uhlherr and Jacques Tits (\cite{Muh-Tit:06}). The $\mathsf{F}_4$ case has been announced by Chris Parker and Katrin Tent (\cite{Par-Ten:*}). Bernhard Leeb and Carlos Ramos Cuevas gave an alternative proof for the $\mathsf{F}_4$ case and also proved the $\mathsf{E}_6$ case (\cite{Lee-Ram:*}). Finally Carlos Ramos Cuevas proved the $\mathsf{E}_7$ case and $\mathsf{E}_8$ case in \cite{Ram:*}. The reducible cases obey the conjecture if their irreducible components do. So all thick spherical buildings obey the center conjecture.

For weak spherical buildings a direct factor of type $\mathsf{H}_3$ or $\mathsf{H}_4$ is possible as well. The center conjecture for the first case follows from results in~\cite{Bal-Lyt:05}, the second is still open.

\section{Main results}
The first main result is a fixed point theorem also valid for $\R$--buildings which are not metrically complete.

\begin{Theorem}\label{main}
A finitely generated bounded group $G$ of isometries of an $\R$--building $\Lambda$ admits a fixed point.
\end{Theorem}

The second main result characterizes metrically complete $\R$--buildings in terms of their wall trees.

\begin{Theorem}\label{main2}
An $\R$--building is metrically complete if and only if its wall trees are metrically complete.
\end{Theorem}

\section{Useful lemmas}
Let $(\Lambda,\cF)$ be an $\R$--building.

\begin{lemma}\label{str}
Let $C_\beta$ and $C_\gamma$ be two sectors based at respectively $\beta$ and $\gamma$, and having the same direction $C_\infty$. Then there exists a constant $k\in \R^+$ depending only on the type of the $\R$--building, such that there exists a point $\delta$ for which the sector $C_\delta$ is a subsector of both $C_\beta$ and $C_\gamma$, and $\dd(\beta,\delta), \dd(\gamma,\delta) \leq k \dd(\beta,\gamma)$.
\end{lemma}
\proof
Embed the sector $C_\beta$ in an apartment $\Sigma$ and the sector $C_\gamma$ in an apartment $\Sigma'$. Let $\delta$ be the point of $C_\beta \cap C_\gamma$ closest to $\beta$ (possible because this intersection is a closed subset of $\Sigma$ due to Condition (A2), and non-empty because parallel sectors  are asymptotic and therefore there exists a subsector contained in both sectors). The sector $C_\delta$ is a subsector of both $C_\beta$ and $C_\gamma$.

Let $D_\infty$ and $D'_\infty$ be the chambers opposite $C_\infty$ in respectively $\Sigma_\infty$ and $\Sigma'_\infty$. Note that $\beta\in D_\delta$ and $\gamma \in D'_\delta$. Due to the way we defined $\delta$, we have that $D_\delta \cap D'_\delta = \{\delta\}$. Consider the retraction $r$ on the apartment $\Sigma$ centered at the germ of $D_\delta$ (see Lemma~\ref{p17}). This retraction maps the sector $D'_\delta$ to some sector $D''_\delta$ in $\Sigma$, only sharing its base point $\delta$ with the sector $D_\delta$. As $r(\gamma)$ lies in $D''_\delta$, it follows that there exists some constant $k$ depending on the minimal angle of halflines in two sectors not sharing simplices, such that $\dd(\beta,\delta), \dd(r(\gamma),\delta) \leq k \dd(\beta,r(\gamma))$. Because the retraction preserves distances to $\delta$, and does not increase the other distances, this implies the desired result.
\qed

\begin{cor}\label{cor:connie}
There exists a constant $k'$ depending only on the type of the $\R$--building, such that for each sector $C_\beta$ and $l \in \R^+$, there exists a point $\delta \in C_\beta$ with $\dd(\beta,\delta)= k'l$, such that for each point $\gamma$ at distance at most $l$ from $\beta$, the sector $C_\delta$ is a subsector of $C_\gamma$.
\end{cor}
\proof
All the sectors $C_\rho$ with $\dd(\rho,\beta)<t$, $t\in \R^+$ and $\rho \in C_\beta$, contain a common point $\tau$ which lies at a distance $k''t$ from $\beta$, with $k''$ some constant. The result then follows from applying the above lemma.
\qed

\begin{cor}\label{cor:equiv}
Let $C_\infty$ and $D_\infty$ be two adjacent chambers at infinity. If for a point $\beta \in \Lambda$ the germs of the sectors $C_\beta$ and $D_\beta$ based at $\beta$ are the same, then there exists an $l >0$ such that for each point $\gamma \in \Lambda$ with $\dd(\beta,\gamma) < l$, the germs of the sectors $C_\gamma$ and $D_\gamma$ based at $\gamma$ are the same.
\end{cor}
\proof
Because the germs of sectors $C_\beta$ and $D_\beta$ are the same, there is an $l' >0$ such that the intersections of the closed ball with radius $l'$ centered at $\beta$ with either the sector $C_\beta$ or $D_\beta$ are the same. This implies that for each point $\delta$ in the intersection of $C_\beta$ and the open ball with radius $l'$ centered at radius $\beta$, the germs of sectors $C_\delta$ and $D_\delta$ are the same. Lemma~\ref{str} now implies that there exists an $l>0$ such that for each point $\gamma \in \Lambda$ with $\dd(\beta,\gamma) < l$, the sectors $C_\beta$ and $C_\gamma$ have a point $\epsilon$ in common at distance strictly less than $l'$ from $\beta$. If the germs of the sectors $C_\gamma$ and $D_\gamma$ were different, then the germs of sectors $C_\epsilon$ and $D_\epsilon$ would be different as well, which is a contradiction. \qed

%Let $A_\infty$ be the common panel of the chambers $C_\infty$ and $D_\infty$.

% are adjacent chambers and the germ of sectors $C_\beta$ and $D_\beta$ are the same

%\qed

\begin{lemma}\label{embed}
The $\R$--building $(\Lambda,\cF)$ can be embedded in an $\R$--building $(\Lambda',\cF')$ of the same type which is metrically complete. Also the action of the isometries of $(\Lambda,\cF)$ can be extended to isometries of $(\Lambda',\cF')$.
\end{lemma}
\proof
This is a (direct) consequence of~\cite[Theorem 5.1.1]{Kle-Lee:97} by taking the $\omega$--limit with respect to the constant scaling
sequence. The resulting space is a complete $\R$--building. (Actually, Kleiner and Leeb assume completeness in their paper, but in this theorem it is irrelevant whether the $\R$--building you start with is complete or not). The choice of the base point doesn't matter and the construction is functorial, so isometries extend. 
\qed

From now on assume that one has an $\R$--building $(\Lambda',\cF')$ as described in the above lemma. Let $\overline{\Lambda}$ be the closure of $\Lambda$ in $\Lambda'$. The metric space defined on $\Lambda$ is complete if and only if $\Lambda = \overline{\Lambda}$. 

Choose $\alpha$ to be a point in $\overline{\Lambda} \backslash \Lambda$ (which is only possible if the metric space defined on $\Lambda$ is not complete). Let $K$ be the chamber subcomplex of the residue $\Lambda'_\alpha$ in $\alpha$ consisting of the germs of sectors-faces based at $\alpha$ with at infinity a simplex of $\Lambda_\infty$. So $K$ is a chamber subcomplex of $\Lambda'_\alpha$. 

\begin{cor}\label{inter}
Let $C_\infty$ be some chamber of the building at infinity $\Lambda_\infty$. Then the interior of the sector $C_\alpha$ lies in $\Lambda$.
\end{cor}
\proof
This follows from applying Corollary~\ref{cor:connie} to a sequence of points in $\Lambda$ converging to $\alpha$.
\qed

\begin{lemma}\label{finder}
Let $C_\alpha$ and $C'_\alpha$ be two sectors based at $\alpha$ with their directions $C_\infty, C'_\infty$ in $\Lambda_\infty$. Let $w$ be the Weyl distance from the germ of $C_\alpha$ to $C'_\alpha$. Then there exists a sector $C''_\alpha$ based at $\alpha$ with the same germ as $C'_\alpha$ such that the direction $C''_\infty$ is a chamber of $\Lambda_\infty$ and that the Weyl distance from $C_\infty$ to $C''_\infty$ is $w$.
\end{lemma}
\proof
Using Lemma~\ref{p8} we know that there exists an apartment $\Sigma$ of $\Lambda'$ containing both the sector $C_\alpha$ and the germ of the sector $C'_\alpha$. Let $\beta$ be a point of $\Sigma$ in the interior of $C'_\alpha$ such that the germ of $C'_\beta$ also lies in $\Sigma$. The Weyl distance from the germ of $C_\beta$ to the germ of $C'_\beta$ is also $w$. Because both the point $\beta$ and the sectors $C_\beta$ and $C'_\beta$ lie in $\Lambda$, one can use Lemma~\ref{p8} again to find an apartment $\Sigma'$ of $\Lambda$ containing $C_\beta$ and the germ of $C'_\beta$. Let $C''_\beta$ be the sector in this apartment $\Sigma'$, based at $\beta$, and with the same germ as $C'_\beta$.  The Weyl distance between $C_\infty$ and $C''_\infty$ is $w$ because the sectors $C_\beta$ and $C''_\beta$ lie in one apartment. As the germs of the sectors $C'_\alpha$ and $C''_\alpha$ are the same, we have proven the lemma.
\qed

\begin{lemma}
The chamber subcomplex $K$ of $\Lambda'_\alpha$ is convex.
\end{lemma}
\proof
Let $C_\alpha$ and $C'_\alpha$ be two sectors with chambers of $\Lambda_\infty$ at infinity. Let the Weyl distance from the germ of $C_\alpha$ to the germ of $C'_\alpha$ be $w$. By the previous lemma, one can assume that the Weyl distance between $C_\infty$ and $C'_\infty$ also is $w$.

Assume we have a germ of sector $D_\alpha$ in the convex hull of the germs of the sectors $C_\alpha$ and $C'_\alpha$. So if the Weyl distance from $C_\alpha$ to $D_\alpha$ is $v$, and the Weyl distance from $D_\alpha$ to $C'_\alpha$ is $v'$ then $w=vv'$. These distances define $D_\alpha$ uniquely. Because the Weyl distance from $C_\infty$ to $C'_\infty$ also is $w$, one can find a chamber $E_\infty$ at infinity such that the Weyl distance from $C_\infty$ to $E_\infty$ is $v$ and the Weyl distance from $E_\infty$ to $C'_\infty$ is $v'$. The distances between the germs $C_\alpha$, $E_\alpha$ and $C'_\alpha$ have to be smaller or equal than these (`smaller' with respect to the word metric on $\overline{W}$ defined by the generating set $S$), but as the distance between $C_\alpha$ and $C'_\alpha$ stays the same, they all stay the same. By uniqueness we can conclude that the germs $[D]_\alpha$ and $[E]_\alpha$ are the same and that $K$ is convex. \qed

\begin{cor}\label{oppie}
No two germs of sectors in $K$ are opposite.
\end{cor}
\proof
Assume that two germs of sectors $C_\alpha$ and $C'_\alpha$ based at $\alpha$ are opposite with $C_\infty,C'_\infty$ chambers of $\Lambda_\infty$. Lemma~\ref{p12} implies that there is an unique apartment $\Sigma$ containing both sectors $C_\alpha$ and $C'_\alpha$. The apartment at infinity will be the unique apartment $\Sigma_\infty$ containing the opposite chambers $C_\infty$ and $C'_\infty$. A consequence is that $\Sigma_\infty$ is an apartment of $\Lambda_\infty$, and so also that $\Sigma$ is an apartment of $\Lambda$. This yields that $\alpha$, being a point of the apartment $\Sigma$, lies in $\Lambda$, which is a contradiction.
\qed

\section{Proof of the first main result}
Let $G$ be a finitely generated bounded group of isometries of an $\R$--building $(\Lambda,\cF)$. Assume that the Main Result~\ref{main} does not hold, or equivalently, that $G$ does not fix points in $\Lambda$ (the Bruhat--Tits fixed point theorem~\ref{bt} implies that $\Lambda$ is not complete if this is the case).

Embed $\Lambda$ in a complete building $\Lambda'$ as described in Lemma~\ref{embed}. The closure $\overline{\Lambda}$ in $\Lambda'$ is a complete CAT(0)--space, so we can apply the Bruhat--Tits fixed point theorem~\ref{bt} and obtain a point $\alpha \in \overline{\Lambda} \subset \Lambda'$ fixed by $G$. Then again as in the previous section, we obtain a convex chamber subcomplex $K$ of the residue $\Lambda'_\alpha$. As $\alpha$ is fixed and $\Lambda$ stabilized by $G$, this group also acts on the (weak) spherical building $\Lambda'_\alpha$ and stabilizes the convex chamber complex $K$. One can now apply the center conjecture~\ref{ccj}, and as Corollary~\ref{oppie} eliminates one option, we know that $G$ stabilizes some non-trivial simplex of $K$. Let $A$ be a maximal stabilized simplex of $K$.

The next step is to investigate the residue of the simplex $A$, which is again a (weak) spherical building, in the spherical building $\Lambda'_\alpha$. The germs in $K$ which contain $A$ form a convex chamber complex of this residue, so we can again apply the center conjecture~\ref{ccj} on this new convex chamber complex.  However, a stabilized non-trivial simplex in the residue is impossible due to the maximality of $A$. So there exist two chambers $C$ and $D$ in $K$, both containing the fixed simplex $A$, and such that the corresponding chambers in the residue of $A$ are opposite. 

Using Lemma~\ref{finder}, one can find two sectors $S_\alpha$ and $S'_\alpha$ lying in one apartment such that their germs equal respectively $C$ and $D$. The intersection $S_\alpha \cap S'_\alpha$ is a sector-face $R_\alpha$ with germ $A$. The interior of both sectors lies in $\Lambda$ due to Corollary~\ref{inter}. Because $\Lambda$ is convex within $\Lambda'$, it follows that the points of $R_\alpha$ not lying on a non-maximal face of this sector-face lie in $\Lambda$.  

Let $L$ be the barycentric closed halfline (with endpoint $\alpha$) of the sector-face $R_\alpha$. From the above discussion it follows that this halfline, except for the point $\alpha$, lies in $\Lambda$. Parametrize this line by a map $\phi:\R^+ \rightarrow L$ such that for each $l \in \R^+$ one has that $\dd(\phi(l),\alpha)=l$. Note that $\phi(0)=\alpha$. As the group $G$ stabilizes the germ $A$, there exists for each element $g \in G$ a positive number $l_g>0$ such that each point $\phi(l)$ with $l \in [0,l_g]$ is fixed by $g$. Let $\{g_1, \dots, g_k\}$ be a finite generating set of $G$ and $l'>0$ be the minimum of the $l_{g_i}$ with $i \in \{1, \dots, k\}$. The point $\phi(l')$, which lies in $\Lambda$, is fixed by a generating set of $G$, and hence by the entire group $G$. So we have proven that there does exist a point in $\Lambda$ fixed by $G$.

\begin{rem} \em
If the Coxeter system has a direct factor of type $\mathsf{H}_4$, then the center conjecture~\ref{ccj} has not been proven yet. However this does not pose a problem for our purposes. If such a case occurs one can restrict the Weyl group of both the spherical building at infinity and the residue at $\alpha$ (which is a retract of the building at infinity) to no longer have a direct factor of type $\mathsf{H}_4$. The convex chamber complex $K$ stays a chamber complex after this restriction of the Weyl group, so we can apply the center conjecture in this case (see~\cite[p. 3]{Ram:*} for more information on this). 
\end{rem}

\section{Proof of the second main result}

First assume that the metric space defined on $\Lambda$ is complete, and let $m$ be a wall of the spherical building at infinity. Let $(\alpha_n)_{n\in \N}$ be a Cauchy sequence in the wall tree $T(m)$. The union of the apartments of the $\R$--building which at infinity contain $m$ forms a subset $K \subset \Lambda$ isometric to the direct product of the metric space formed by $T(m)$ and $\R$.

Using this subset $K$, we can `lift' the Cauchy sequence $(\alpha_n)_{n\in \N}$ to a Cauchy sequence $(\beta_n)_{n\in \N}$ in $K \subset \Lambda$. As the metric space defined on $\Lambda$ is complete, this sequence converges to some point $\beta \in \Lambda$. Our goal is to prove that the point $\beta$ lies in $K$, implying that the original sequence $(\alpha_n)_{n\in \N}$ converges. For this we have to prove that $\beta$ lies in an apartment which at infinity contains the wall $m$. Let $S_\infty$ and $S'_\infty$ be two opposite maximal sector-panels of $m$; if we can prove that the germs of sector-panels $[S]_\beta$ and $[S']_\beta$ in the residue at $\beta$ are still opposite, we are done (because Lemma~\ref{p12} then implies that $\beta$ lies in an apartment containing $m$ at infinity). Equivalent with this last statement is that for a shortest gallery from a chamber $C_\infty$ containing $S_\infty$ to a chamber $C'_\infty$ containing $S'_\infty$ (`shortest' meaning minimal over all choices of $C_\infty$ and $C'_\infty$), the corresponding gallery in the residue $\Lambda_\beta$ between the germs of sectors $[C]_\beta$ to $[C']_\beta$ always is non-stammering. As this is the case for each point of $K$, and hence each point of the sequence $(\beta_n)_{n\in \N}$, Corollary~\ref{cor:equiv} implies that this is also the case for $\beta$. So we have proven that the metric space defined by the $\R$-tree $T(m)$ is complete.

We are now left with the other direction to prove. Assume that all the trees corresponding to walls at infinity are complete. Let $ (\alpha_n)_{n\in \N}$ be a Cauchy sequence in the metric space $(\Lambda,\dd)$. Using Lemma~\ref{embed}, we embed $(\Lambda,\cF)$ in a complete $\R$--building $(\Lambda',\cF')$. The Cauchy sequence $ (\alpha_n)_{n\in \N}$ then converges to some point $\alpha \in \overline{\Lambda}$. We need to prove that $\alpha$ is a point of $\Lambda$, so assume this is not the case. 

%Choose some chamber $C_\infty$ at infinity and consider the sequence of sectors $(C_{\alpha_n})_{n\in \N}$. Using Corollary~\ref{cor:connie} one can find a sequence of points $(\beta_n)_{n\in \N}$ in $\Lambda$ which also converges to the point $\alpha$, and such that if $i<j$, then the sector $C_{\beta_i}$ is a subsector of $C_{\beta_j}$. In the completion $\bar{\Lambda}$ we obtain a subset isometric to a sector, where the `source' is $\alpha$ (by applying Corollaire 2.11 from~\cite{Par:00} and its preceding text). Note that the interior (as it would be in an apartment) lies in $\Lambda$.

Choose a sector $C_\alpha$ based at $\alpha$. Due to Corollary~\ref{inter} the interior of the sector $C_\alpha$ lies in $\Lambda$. One can find a sequence of points $(\beta_n)_{n\in \N}$ in the interior of the sector $C_\alpha$ which also converges to the point $\alpha$.

Let $P_\infty$ be a panel of $C_\infty$. The sequence $(P_{\beta_n})_{n\in \N}$ of parallel sector-panels forms a Cauchy sequence in the panel tree $T(P_\infty)$, contained in an open halfline. Using the completeness of this tree, we can embed this open halfline into a closed halfline, and then into to an apartment (essentially using Lemma~\ref{p8}), and find a chamber $C'_\infty$ in $\Lambda_\infty$ containing $P_\infty$ such that the germs of the sectors $C_{\beta_n}$ and $C'_{\beta_n}$ are not the same for all $n \in \N$. It follows that the germs of the sectors $C_\alpha$ are $C'_\alpha$ are not the same, but adjacent, having the germ of the sector-panel $P_\alpha$ in common.

Repeating this algorithm one can obtain two sectors based at $\alpha$ with at infinity chambers of $\Lambda_\infty$ and opposite germs, but this is in contradiction with Corollary~\ref{oppie} and $\alpha \notin \Lambda$. This proves the second main result.

%%%%%%%%%%%%%%%%%%%%   End of main body of article
%
%                             References
%
%   BiBTeX users uncomment the following line:
%
%\bibliographystyle{gtart}
%


\begin{thebibliography}{99}
\bibitem{Abr-Bro:08}
P. Abramenko and K. Brown, {\em Buildings: Theory and applications}, Springer, 2008.

\bibitem{Bal-Lyt:05}
A. Balser and A. Lytchak, Centers of convex subsets of buildings, \emph{Ann. Glob. Anal. Geom.} \textbf{28}, No. 2 (2005), 201--209

\bibitem{Ber-Kap:*} 
A.~Berenstein and M.~Kapovich, Affine buildings for dihedral groups, \emph{preprint}.


\bibitem{Bru-Tit:72} 
F. Bruhat and J. Tits, Groupes r\'eductifs sur un corps local,
I.~Donn\'ees radicielles valu\'ees, {\em Inst.\ Hautes
\'Etudes Sci.\ Publ.\ Math}.\ {\bf 41} (1972), 5--252.

\bibitem{Par-Ten:*}
C. Parker and K. Tent, Convexity in buildings, talk at the Oberwolfach conference \emph{Buildings: Interactions with algebra and geometry}, Oberwolfach report no. 3/2008, 33--34.
%C. Parker and K. Tent, Completely reducible subcomplexes of spherical buildings, preprint.

\bibitem{Par:00}
A. Parreau, Immeubles affines: construction par les normes et \'etude des isom\'etries, in \emph{Crystallographic groups
and their generalizations (Kortrijk, 1999)}, Contemp. Math. \textbf{262}, Amer. Math. Soc., Providence, RI,
2000, pp. 263--302.

\bibitem{Kle-Lee:97}
B. Kleiner and B. Leeb, Rigidity of quasi-isometries for symmetric spaces and Euclidean buildings, \emph{Inst. Hautes \'Etudes Sci. Publ. Math.}, \textbf{86} (1997), 115--197.

\bibitem{Lee-Ram:*}
B. Leeb and C. Ramos Cuevas, The center conjecture for spherical buildings of types $\mathsf{F}_4$ and $\mathsf{E_6}$, \emph{preprint}.

\bibitem{Muh-Tit:06}
B. M\"uhlherr and J. Tits, The center conjecture for non-exceptional buildings, \emph{Journal of Algebra}, \textbf{300} (2006), 687--706.

\bibitem{Ram:*}
C. Ramos Cuevas, The center conjecture for thick spherical buildings, \emph{preprint}.

\bibitem{Str-Mal:*}
K. Struyve and H. Van Maldeghem, Two-dimensional affine $\R$--buildings defined by generalized polygons with non-discrete valuation, \emph{Pure Appl. Math Q.}, to appear.

\bibitem{Tit:86}
J. Tits, Immeubles de type affine, In \emph{Buildings and the Geometry of Diagrams, Springer Lecture Notes}
\textbf{1181} (Rosati ed.), Springer Verlag, 1986, pp. 159--190.


\end{thebibliography}
\end{document}